\documentclass[a4paper]{article}
\usepackage{amsfonts,amssymb}
\usepackage[matrix,arrow,curve]{xy}

\newtheorem{dfn}{Definition}[section]
\newtheorem{prop}[dfn]{Proposition}
\newtheorem{theo}[dfn]{Theorem}
\newtheorem{ex}[dfn]{Example}
\newtheorem{exs}[dfn]{Examples}
\newtheorem{cor}[dfn]{Corollary}
\newtheorem{lem}[dfn]{Lemma}

\newcommand{\Rem}{\noindent {\sc Remark. }}

\newcommand{\Proof}{\noindent {\sc Proof. }}
\newcommand{\eop}{$\square$}

\newcommand{\ra}{\rightarrow}
\newcommand{\lra}{\longrightarrow}
\newcommand{\oo}{\,\mbox{-}\,}
\newcommand{\com}{\circ}

\newcommand{\RR}{\mathbb{R}}

\newcommand{\cA}{\mathord{\mathcal{A}}}
\newcommand{\cG}{\mathord{\mathcal{G}}}
\newcommand{\cH}{\mathord{\mathcal{H}}}
\newcommand{\cV}{\mathord{\mathcal{V}}}
\newcommand{\cL}{\mathord{\mathcal{L}}}
\newcommand{\cF}{\mathord{\mathcal{F}}}
\newcommand{\cK}{\mathord{\mathcal{K}}}

\newcommand{\fX}{\mathord{\mathfrak{X}}}
\newcommand{\fg}{\mathord{\mathfrak{g}}}
\newcommand{\fh}{\mathord{\mathfrak{h}}}
\newcommand{\fk}{\mathord{\mathfrak{k}}}

\newcommand{\id}{\mathord{\mathit{id}}}
\newcommand{\pr}{\mathord{\mathit{pr}}}

\newcommand{\eC}{\mathord{\mathit{C}^{\infty}}}

\newcommand{\src}{\mathord{\alpha}}
\newcommand{\trg}{\mathord{\beta}}
\newcommand{\mlt}{\mathord{\mu}}
\newcommand{\inv}{\mathord{\mathit{inv}}}
\newcommand{\uni}{\mathord{\mathit{uni}}}
\newcommand{\anchor}{\mathord{\rho}}
\newcommand{\act}{\mathord{\vartheta}}

\newcommand{\rank}{\mathord{\mathrm{rank}}}
\newcommand{\Der}{\mathord{\mathrm{Der}}}

\newcommand{\Ker}{\mathord{\mathrm{Ker}}}
\newcommand{\Mon}{\mathord{\mathrm{Mon}}}
\newcommand{\Hol}{\mathord{\mathrm{Hol}}}
\newcommand{\Iso}{\mathord{\mathrm{Iso}}}
\newcommand{\Ad}{\mathord{\mathrm{Ad}}}

\title{\bf On integrability of infinitesimal actions}
\author{{\sc Ieke Moerdijk} \\
  {\normalsize Mathematical Institute, Utrecht University}
   \vspace{-1mm} \\
  {\normalsize P.O. Box 80.010, 3508 TA Utrecht,
               The Netherlands} \\
  {\small \it Email address: \tt moerdijk@math.uu.nl}
              \vspace{3mm} \\
  {\sc Janez Mr\v{c}un} \\
  {\normalsize Department of Mathematics, University of Ljubljana}
               \vspace{-1mm} \\
  {\normalsize Jadranska 21, 1000 Ljubljana, Slovenia} \\
  {\small \it Email address: \tt Janez.Mrcun@fmf.uni-lj.si}}
\date{}

\begin{document}

\maketitle

\begin{abstract}
We use foliations and connections on principal Lie groupoid
bundles to prove various integrability results for Lie algebroids.
In particular, we show, under quite general assumptions, that the
semi-direct product associated to an infinitesimal action of one
integrable Lie algebroid on another is integrable. This generalizes
recent results of Dazord and Nistor.
\end{abstract}

\section*{Introduction} \label{sec:introduction}

Lie algebroids have recently turned out to be very useful in several
ways, e.g. related to deformation quantization of manifolds and to
Poisson geometry (see for example \cite{CdSW,CF,K,L,NT}).
The notion of a Lie algebroid itself goes back to Pradines
\cite{Pr,Pr1}. He constructed for each Lie groupoid such a Lie
algebroid, and outlined a Lie theory describing a correspondence
between Lie groupoids and algebroids
completely analogous to the classical theory for Lie groups and
algebras. It remained a problem to develop the details of this
theory, until Almeida and Molino \cite{AM} provided a counterexample
to one of Pradines main assertions: they proved that a transversely
complete foliation gives rise to a transitive Lie algebroid which is
integrable if and only if the foliation is developable.
Nevertheless, there are some well-known classical examples of Lie
algebroids which can be integrated. For example, Douady and Lazard
\cite{DL} proved that any bundle of Lie algebras can be integrated
to a bundle of Lie groups, which in the language of Lie algebroids
means that any Lie algebroid with trivial anchor map is integrable.

Recently, several positive integrability results have been discovered.
Dazord \cite{D} proved that the
transformation Lie algebroid associated to an infinitesimal action
of a Lie algebra on a manifold is integrable to a Lie groupoid.
Furthermore, using the integrability of foliation algebroids
and the result of Douady and Lazard, Nistor \cite{N}
proved that any regular Lie algebroid which admits a flat splitting
is integrable. His motivation for this result was to construct
examples of pseudodifferential operators on groupoids \cite{NWX}.
In a recent thesis, Debord proved that any Lie algebroid with almost
injective anchor map is integrable \cite{Debord}.

The purpose of this paper is to prove some integrability
results, which include the ones of Dazord and Nistor. More specifically,
our main results concern actions of one Lie algebroid on another.
We prove that in many cases, the semi-direct product of such an
infinitesimal action (described in \cite{HM}) is integrable whenever
each of the algebroids is.

The outline of this paper is as follows. In the first section,
we recall some basic definitions and main examples concerning Lie algebroids
and Lie groupoids, and fix the notations. In the second section,
we discuss principal $G$-bundles for a Lie groupoid $G$, and introduce a notion
of connection which takes values in the Lie algebroid of $G$.
In the third section, we first
give a quick and uniform treatment of some basic results
of the Lie theory for groupoids.
More specifically,
we construct for each Lie groupoid a
source-simply connected cover having the same algebroid
(the source-simply connected Lie groupoids play the role analogous
to the simply connected Lie groups). The existence of
such a cover was proved earlier by more involved methods
in some special cases, e.g. \cite{BM,M}; see \cite{M1} for a survey.
Furthermore, we use foliation theory and connections to show that integrability
is inherited by subalgebroids, and to give
a simple proof of Mackenzie and Xu's result
concerning integrability of morphisms between algebroids \cite{MX}.
In this section, we also prove an integrability theorem stating
that an action of a Lie groupoid $G$ on an
integrable Lie algebroid $\fh$
can be integrated to an action of $G$ on the integral groupoid of $\fh$.

In section 4, we discuss derivations on Lie algebroids.
For the Lie algebroid $\fh$ of a source-simply connected groupoid $H$,
we prove that the Lie algebra $\Der(\fh)$ of derivations on $\fh$
is isomorphic to the Lie algebra of multiplicative vector fields on $H$. We also
discuss actions of another Lie algebroid $\fg$ on $\fh$ in terms
of the algebra $\Der(\fh)$, and recall the construction of semi-direct
products from \cite{HM}.

The main results of this paper are contained in section 5. Here we prove that
under quite general assumptions, the semi-direct product of integrable Lie
algebroids is itself integrable. More specifically, this holds if the
algebroid $\fh$ is a foliation (Theorem \ref{iaf1}), or
if $\fh$ is ``proper'' over $\fg$ in a suitable sense
(Corollary \ref{iapm2}), or if $\fh$ is integrable
by a source-compact source-simply connected Lie groupoid
(Theorem \ref{iascg1}).
At present, we do not know whether in general $\fg\ltimes\fh$ is integrable
if $\fg$ and $\fh$ are.
\vspace{3mm}

\noindent
{\bf Acknowledgements.}
We are grateful to E. van den Ban, J. Duistermaat and
K. Mackenzie for helpful discussions. We acknowledge support of the
Dutch Science Foundation (NWO), which supported a visit of the second author
to the University of Utrecht in the fall of 1999, when a first draft of
this paper was written.

\section{Lie algebroids and Lie groupoids}

\subsection{Lie algebroids} \label{sec:alg}

In this section we recall the definition of a Lie algebroid as well as
some of the main examples. For a more extensive discussion we refer to
\cite{CdSW,M}. Throughout this paper, we shall work in the smooth
context, so ``manifold''
means smooth manifold, ``map'' means smooth map, ``vector bundle''
means smooth real vector bundle, etc.

Let $M$ be a manifold.
A {\em Lie algebroid} over $M$
is a smooth vector bundle $\pi:\fg\ra M$,
together with a map $\anchor:\fg\ra T(M)$
of vector bundles over $M$
and a (real) Lie algebra structure $[\oo,\oo]$ on the
vector space $\Gamma\fg$ of sections of $\fg$ such that
\begin{enumerate}
\item  [(i)] the induced map $\Gamma(\anchor):\Gamma\fg\ra\fX(M)$
             is a Lie algebra homomorphism, and
\item [(ii)] the Leibniz identity
             $$ [X,fX']=f[X,X']+\Gamma(\anchor)(X)(f)X' $$
             holds for any $X,X'\in\Gamma\fg$ and any $f\in\eC(M)$.
\end{enumerate}
The map $\anchor$
is called the {\em anchor} of the Lie algebroid $\fg$.
The map $\Gamma(\anchor)$ is often denoted by $\anchor$ as well, and also called
the anchor. The manifold $M$ is called the {\em base manifold}
of the Lie algebroid $\fg$.

\begin{exs} \rm \label{alg1}
(i) Every finite dimensional
Lie algebra is a Lie algebroid over a one point space.

(ii) Any manifold $M$ can be viewed as a Lie algebroid in two ways, by
taking the zero bundle over $M$ (which we shall denote simply by $M$),
or by taking the tangent bundle over $M$ with
the identity map for the anchor (we shall denote this Lie algebroid by $T(M)$).

(iii) Any vector bundle $E$ over $M$ can be viewed as a Lie algebroid over $M$,
with zero bracket and anchor.

(iv) A vector bundle $E$ over $M$ with a smoothly varying Lie
algebra structure on its fibers (i.e. a bundle of Lie algebras)
can be viewed as a Lie algebroid over $M$ with zero anchor.

(v) A foliation $\cF$ of $M$ is by definition an involutive (hence integrable)
subbundle of $T(M)$. Thus a foliation of $M$ is the same thing as a Lie
algebroid over $M$ with injective anchor map.

(vi) Let $M$ be a manifold equipped with an infinitesimal action of a Lie
algebra $\fg$, i.e. a Lie algebra homomorphism $\gamma:\fg\ra\fX(M)$.
The trivial bundle $\fg\times M$ over $M$ has the structure of a Lie algebroid,
with anchor given by $\anchor(\xi,x)=\gamma(\xi)_{x}$, and Lie bracket
$$ [u,v](x)=[u(x),v(x)]+(\gamma(u(x))(v))(x)-(\gamma(v(x))(u))(x)\;,$$
for $u,v\in\eC(M,\fg)\cong\Gamma(M,\fg\times M)$ and $x\in M$. This Lie algebroid
is called the {\em transformation} algebroid associated to the
infinitesimal action.

(vii) Let $(M,\Pi)$ be a Poisson manifold. Then there is a natural
Lie algebra structure on $\Omega^{1}(M)$ which makes $T^{\ast}(M)$ into a Lie
algebroid over $M$. The anchor of this algebroid is $-\tilde{\Pi}$, where
$\tilde{\Pi}:T^{\ast}(M)\ra T(M)$ is induced by the bivector field $\Pi$.
For details, see e.g. \cite{CdSW}.
\end{exs}

Let $\fg$ be a Lie algebroid over $M$ and
$\phi:N\ra M$ a map of manifolds.
Consider the pull-back bundle $\phi^{\ast}\fg$.
The sections of the form $\phi^{\ast}(X)=(\id,X\com\phi)$, $X\in\Gamma\fg$,
span $\Gamma\phi^{\ast}\fg$ as a $\eC(N)$-module. In fact,
the map $\eC(N)\otimes_{\eC(M)}\Gamma\fg\ra\Gamma\phi^{\ast}\fg$,
which sends $f\otimes X$ to $f\phi^{\ast}(X)$, is an isomorphism
\cite{GHV}.
$$
\xymatrix{
\phi^{\ast}\fg =N\times_{M}\fg \ar[r] \ar@<2pt>[d] & \fg \ar@<2pt>[d]^{\pi} \\
N \ar@<2pt>[u]^{\phi^{\ast}(X)} \ar[r]^{\phi} & M \ar@<2pt>[u]^{X}
}
$$
Let $\fh$ be a Lie algebroid over $N$,
and $\Phi:\fh\ra\fg$ a bundle map over $\phi:N\ra M$.
A {\em $\Phi$-decomposition} of a section $Y\in\Gamma\fh$ is a decomposition
$(\id,\Phi\com Y)=\sum_{i}f_{i}\phi^{\ast}(X_{i})\in\Gamma\phi^{\ast}\fg$,
for some $f_{i}\in\eC(N)$ and $X_{i}\in\Gamma\fg$.
Such a bundle map $\Phi:\fh\ra\fg$ over $\phi$ is a
{\em morphism} of Lie algebroids \cite{HM} if it preserves the anchor,
i.e. $\anchor\com\Phi=d\phi\com\anchor$, and
if it preserves the bracket in the following sense: for any
$Y\in\Gamma\fh$ with a $\Phi$-decomposition
$\sum_{i}f_{i}\phi^{\ast}(X_{i})$, and any $Y'\in\Gamma\fh$
with a $\Phi$-decomposition
$\sum_{j}f'_{j}\phi^{\ast}(X'_{j})$,
$$ \sum_{i,j}f_{i}f'_{j}\phi^{\ast}([X_{i},X'_{j}])
   + \sum_{j}\anchor(Y)(f'_{j})\phi^{\ast}(X'_{j})
   - \sum_{i}\anchor(Y')(f_{i})\phi^{\ast}(X_{i}) $$
is a $\Phi$-decomposition of $[Y,Y']$.

\subsection{Lie groupoids} \label{sec:grp}

Like the previous one, this section only serves to
recall some basic definitions and fix the notations.

A groupoid is a small category $G$ in which all the arrows are invertible.
We shall write $G_{0}$ for the set of objects of $G$, while the set of arrows
of $G$ will be denoted by $G_{1}$. We shall often identify
$G_{0}$ with the subset of units of $G_{1}$. The structure
maps of $G$ will be denoted as follows:
$\src,\trg:G_{1}\ra G_{0}$ will stand for the source (domain) map,
respectively the target (codomain) map,
$\mlt:G_{1}\times_{G_{0}} G_{1}\ra G_{1}$ ($\mlt(g,g')=gg'$) for the multiplication
(composition) map, $\inv:G_{1}\ra G_{1}$ ($\inv(g)=g^{-1}$)
for the inverse map and $\uni:G_{0}\ra G_{1}$ ($\uni(x)=1_{x}$)
for the unit map. We sometimes say that $G$ is a groupoid over $G_{0}$.

A {\em Lie groupoid} is a groupoid $G$, equipped with the
structure of smooth manifold both
on the set of arrows $G_{1}$ and on
the set objects $G_{0}$, such that all the structure
maps of $G$ are smooth and $\src$ is a submersion.
Note that this implies that $\trg$ is a submersion as well, that there
is a natural smooth structure on the domain $G_{1}\times_{G_{0}}G_{1}$ of the
multiplication, and that $\uni$ is an embedding.
We shall assume that the manifold of objects $G_{0}$ and the $\src$-fibers
$\src^{-1}(x)$, $x\in G_{0}$, are Hausdorff, but we do not assume
$G_{1}$ to be Hausdorff (cf. Example \ref{gpd1} (v)).

A morphism of Lie groupoids $F:H\ra G$ is a functor which is smooth
as a map between the manifolds of arrows ($F_{1}:H_{1}\ra G_{1}$)
and as a map between the manifolds of objects ($F_{0}:H_{0}\ra G_{0}$).

\begin{exs} \rm \label{gpd1}
(i) Every Lie group is a Lie groupoid over a one point space.

(ii) Any submersion $M\ra B$ gives a Lie groupoid $M\times_{B}M$ over
$M$, with $\src=\pr_{2}$ and $\trg=\pr_{1}$.
In particular, to any manifold $M$ we associate two Lie groupoids:
the one corresponding to the identity $M\ra M$, and the one (the pair groupoid)
corresponding to the map from $M$ to a space with only one point.

(iii) A Lie groupoid $G$ with $\src=\trg$ is just a family of Lie groups
smoothly parametrized by $G_{0}$. In particular, any vector bundle is
a Lie groupoid with $\src=\trg$.

(iv) The fundamental groupoid of a manifold is a Lie groupoid.

(v) A foliated manifold $(M,\cF)$ gives rise to two Lie groupoids, the
holonomy groupoid $\Hol(M,\cF)$ and the monodromy (or homotopy)
groupoid $\Mon(M,\cF)$. The manifold of objects is $M$ in both cases.
If $x,y\in M$ are on different leaves, there are no arrows from $x$ to
$y$ in both cases. If $x$ and $y$ are on the same leaf $L$, the arrows
from $x$ to $y$ in $\Mon(M,\cF)$ are homotopy classes of paths from $x$ to
$y$ inside $L$. Thus $\Mon(M,\cF)$ is the union of the fundamental groupoids
of the leaves, equipped with a suitable smooth structure. The holonomy groupoid
is a quotient of the monodromy groupoid: for $x,y\in L$,
the arrows from $x$ to $y$ in $\Hol(M,\cF)$ are the holonomy classes of
paths from $x$ to $y$ inside $L$.
These monodromy and holonomy groupoids are generally non-Hausdorff.
For details, see \cite{Ph,W}.

(vi) If $M$ is a manifold equipped with a smooth left action of a Lie group
$G$, the translation groupoid $G\ltimes M$ has $M$ for its manifold of objects
and $G\times M$ for its manifold of arrows, an arrow from $x$ to $y$ being
a pair $(g,x)$ with $gx=y$. The multiplication in $G\ltimes M$ is defined by
$(g,x)(g',x')=(gg',x')$ when $x=g'x'$.

(vii) Let $G$ be a Lie group and let $P\ra M$ be a right principal $G$-bundle
over $M$. The pair groupoid $P\times P$ over $P$ of Example \ref{gpd1} (ii)
has a natural diagonal $G$-action,
and the quotient $(P\times P)/G$ is a Lie groupoid over $M$, called the
{\em gauge} groupoid of the principal bundle.
\end{exs}

\subsection{The Lie algebroid of a Lie groupoid}  \label{sec:lalg}

The construction of a Lie algebra $\fg$ of a given Lie group $G$ extends to
groupoids \cite{Pr,M}. Explicitly, if $G$ is a Lie groupoid,
the vector bundle $T^{\src}(G_{1})=\Ker(d\src)$ over $G_{1}$ of
$\src$-vertical tangent vectors pulls back along
$\uni:G_{0}\ra G_{1}$ to a vector bundle $\fg$ over $G_{0}$.
This vector bundle has the structure of a Lie algebroid.
Its anchor $\anchor:\fg\ra T(G_{0})$ is induced by the differential
of the target map, $d\trg:T(G_{1})\ra T(G_{0})$. To define the bracket,
note that the bundle $T^{\src}(G_{1})$ over $G_{1}$ has a right $G$-action,
and that
the $G$-invariant sections of $T^{\src}(G_{1})$ over $G_{1}$
form a Lie subalgebra of $\fX(G_{1})$, which we denote by
$\fX^{\src}_{inv}(G)$. Now the sections of $\fg$ over $G_{0}$
can be identified with the $G$-invariant sections of $T^{\src}(G_{1})$ over
$G_{1}$. Explicitly, a section $X$ of $\fg$ gives an invariant $\src$-vertical
vector field on $G_{1}$ with value
$$ X_{\trg(g)}g $$
at an arrow $g\in G_{1}$.
This identification gives us a Lie algebra structure on
$\Gamma\fg$; see \cite{M} for details. We denote the Lie algebroids
associated to $G$, $H$, etc. by $\fg$, $\fh$, etc, or sometimes by
$\cA(G)$, $\cA(H)$, etc.
The differential of a morphism $F:H\ra G$ of Lie groupoids induces a morphism
$\cA(F):\cA(H)\ra \cA(G)$ of Lie algebroids over $F_{0}:H_{0}\ra G_{0}$,
in a functorial way \cite{HM}.

\begin{dfn} \label{lalg1}
A Lie algebroid $\fg$ is called integrable if it is isomorphic
to the Lie algebroid associated to a Lie groupoid $G$.
If this is the case, then $G$ is called an integral of $\fg$.
\end{dfn}

\Rem
Contrary to the case of finite dimensional Lie algebras and Lie groups,
there exist Lie algebroids which are not integrable. See
\cite{AM,Mo} for an example.

\begin{exs} \rm \label{lalg2}
(i) Any finite dimensional Lie algebra is an integrable Lie algebroid
(``Lie's third theorem'').

(ii) Any foliation $\cF$ of a manifold $M$ is an integrable Lie algebroid,
because we have $\cF\cong\cA(\Hol(M,\cF))\cong\cA(\Mon(M,\cF))$.

(iii) Any bundle $E$ of Lie algebras over $M$ is an integrable Lie
algebroid \cite{DL}. In fact, there exists a bundle of Lie groups
(which may not be locally trivial nor Hausdorff)
which integrates $E$.
\end{exs}

\section{Groupoid bundles and connections}

\subsection{Actions by Lie groupoids}  \label{sec:ag}

In this section we begin by recalling several standard notions and terminology
concerning groupoid bundles which play a central role in the proofs of
some of our integrability results.

Let $G$ be a fixed Lie groupoid, $N$ a manifold and $\epsilon:N\ra G_{0}$
a map. A left {\em action} of $G$ on $N$ {\em along} $\epsilon$
is given by a map $\act:G_{1}\times_{G_{0}}N\ra N$ (we write $\act(g,y)=gy$),
defined on the pull-back
$G_{1}\times_{G_{0}}N=\{(g,y)\,|\,\src(g)=\epsilon(y)\}$, which
satisfies the following identities:
$\epsilon(gy)=\trg(g)$, $1_{\epsilon(y)}y=y$ and
$g'(gy)=(g'g)y$, for any  $g',g\in G_{1}$ and $y\in N$ with
$\src(g')=\trg(g)$ and $\src(g)=\epsilon(y)$.
For such an action one can form the {\em translation} groupoid
$G\ltimes N$ over $N$, with $(G\ltimes N)_{1}=G_{1}\times_{G_{0}}N$,
exactly as for groups (Example \ref{gpd1} (vi)).
We define the quotient $N/G$ as the space of orbits of the groupoid
$G\ltimes N$. This space is in general not a manifold.
A right action of $G$ on $N$ is defined analogously.

Suppose that we have a right $G$-action on a manifold $P$.
If $P$ is equipped with a map $\pi:P\ra B$ and the action is fiberwise in
the sense that $\pi(pg)=\pi(p)$ whenever $pg$ is defined, then $P$ is called
a $G$-bundle over $B$. This $G$-bundle is said to be {\em principal}
if $\pi$ is a surjective submersion and the map
$(\act,\pr_{1}):P\times_{G_{0}}G_{1}\ra P\times_{B}P$ is a diffeomorphism.
In this case the translation groupoid $P\rtimes G$ is isomorphic
to the pair groupoid $P\times_{B} P$ over $P$.
Note that $P/G\cong B$, so $P/G$ is a manifold.

There is also a notion of a (left) action of a Lie groupoid $G$ on
another Lie groupoid $H$. It is given by two (left) actions of
$G$ on $H_{1}$ and on $H_{0}$,
such that the groupoid structure maps of $H$ are
compatible with the actions by $G$ \cite{HM}.
If we denote the action maps on $H_{i}$ by
$\epsilon_{i}:H_{i}\ra G_{0}$ and $\act_{i}:G_{1}\times_{G_{0}}H_{i}\ra H_{i}$,
$i=0,1$, this implies in particular that
$\epsilon_{0}\com\src=\epsilon_{1}=\epsilon_{0}\com\trg$.
Thus the fibers $H_{x}=\epsilon_{1}^{-1}(x)$ are full subgroupoids
of $H$ over $\epsilon_{0}^{-1}(x)$, $x\in G_{0}$.
These are Lie groupoids if $\epsilon_{0}$ is a submersion, and
for each arrow $g\in G_{1}(x',x)$ the action provides
an isomorphism $H_{x'}\ra H_{x}$ of Lie groupoids.

For such an action of $G$ on $H$, one can form the
{\em semi-direct product} groupoid
$$ G\ltimes H $$
over $H_{0}$. For $y,y'\in H_{0}$, an arrow from $y'$ to $y$ in
$G\ltimes H$ is a pair $(g,h)$, where $g$ is an arrow in
$G(\epsilon_{0}(y'),\epsilon_{0}(y))$ and $h$ is an arrow
in $H(gy',y)\subset H_{\epsilon_{0}(y)}$. These arrows compose by
the usual formula
$$ (g,h)(g',h')=(gg',h(gh'))\;.$$
The groupoid $G\ltimes H$ has the natural structure of a Lie groupoid,
as one sees, e.g. when the space of arrows is considered as the fibered product
$$ H_{1}\times_{G_{0}}G_{1}=\{(h,g)\,|\,\epsilon_{0}(\trg(h))=\trg(g)\}\;.$$

\begin{lem}  \label{ag1}
Consider an action of a Lie groupoid $G$ on a Lie groupoid $H$. If $H_{0}$
is a principal $G$-bundle over $B$, then $H_{1}/G$ is a Lie groupoid
over $B\cong H_{0}/G$.
\end{lem}
\Proof
The only thing that has to be shown is that $H_{1}/G$ is a manifold.
We can specify the manifold structure locally in $B$, so it suffices
to consider the case where $\pi:H_{0}\ra B$ has a section $s$. But then
$H_{1}/G$ is isomorphic to the pull-back of $\src:H_{1}\ra H_{0}$ along
$s:B\ra H_{0}$, hence is a manifold. Moreover, this manifold structure is
independent of the choice of $s$, since by principality of the action
on $H_{0}$ any two sections $s$ and $s'$ are related by a map
$\theta:B\ra G_{1}$ as $\theta(b)s(b)=s'(b)$ for all $b\in B$.
Then the same multiplication by $\theta$ establishes a diffeomorphism
between the pull-back of $\src$ along $s$ and the one along $s'$.
\eop
\vspace{3mm}

\Rem
We denote the Lie groupoid $H_{1}/G$ over $B$ by $H/G$.
The quotient morphism $H\ra H/G$ induces for each $x\in H_{0}$ an isomorphism
of $\src$-fibers $\src^{-1}(x)\ra\src^{-1}(\pi(x))$. More precisely,
the square
$$
\xymatrix{
H_{1} \ar[d]_{\src} \ar[r] & H_{1}/G \ar[d]^{\src} \\
H_{0} \ar[r]^-{\pi} & B=H_{0}/G
}
$$
is a pull-back of smooth manifolds.

\subsection{Connections on principal groupoid bundles} \label{sec:cpb}

Let $G$ be a Lie groupoid with Lie algebroid $\fg$,
and let $\pi:P\ra B$ be a principal $G$-bundle along $\epsilon:P\ra G_{0}$.
For any $p\in P$, denote by $\cV_{p}$ the space $\Ker((d\pi)_{p})$
of vertical tangent vectors at $p\in P$.
Thus $\cV$ is an integrable subbundle of $T(P)$. The diffeomorphism
$L_{p}:G(\epsilon(p),\oo)\ra P_{\pi(p)}$, given by $L_{p}(g)=pg^{-1}$, induces
and isomorphism $dL_{p}:\fg_{\epsilon(p)}\ra \cV_{p}$.

Recall that a {\em local bisection} of $G$ is a local section
$s:U\ra G_{1}$ of the source map $\src$, defined on an open
subset $U$ of $G_{0}$,
such that $\trg\com s$ is an open embedding.
We say that $s$ is a local bisection of $G$ through $g\in G_{1}$ if
$g\in s(U)$. It is easy to see that there exist local bisections
through any arrow of $G$. Such a bisection $s$
induces a diffeomorphism
$R_{s}:\epsilon^{-1}(\trg(s(U))\ra\epsilon^{-1}(U)$
by $R_{s}(p)=ps(x)$,
where $x\in G_{0}$ is the unique point with
$\trg(s(x))=\epsilon(p)$.
If $\xi$ is in the kernel of $(d\epsilon)_{p}$, then $(dR_{s})_{p}(\xi)$
depends only on the value $g$ of $s$ with $\trg(g)=\epsilon(p)$, and we will
write $\xi g=(dR_{s})_{p}(\xi)$. Indeed, in this case
one may define $\xi g$ as the image of $(\xi,0)$ along the derivative
of the action $P\times_{G_{0}}G_{1}\ra P$ at $(p,g)$.

Let $\cF$ be a foliation of $B$. Then $\pi^{\ast}(\cF)$ is a foliation
of $P$. An {\em $\cF$-partial connection} on $P$ is a subbundle $\cH$ of
$\pi^{\ast}(\cF)\subset T(P)$ which satisfies the following conditions:
\begin{enumerate}
\item [(i)]   $\pi^{\ast}(\cF)=\cV\oplus \cH$,
\item [(ii)]  $(d\epsilon)(\cH)=0$, and
\item [(iii)] $\cH_{pg}=\cH_{p} g$ for any $p\in P$ and $g\in G$ with
              $\epsilon(p)=\trg(g)$.
\end{enumerate}
Note that $\cH_{p} g$ is well-defined precisely because of the condition (ii)
above. The connection $\cH$ is called {\em flat} if it is integrable.
With a given connection $\cH$, any tangent vector $\xi\in \pi^{\ast}(\cF)_{p}$
has a unique decomposition as a sum $\xi=\xi^{v}+\xi^{h}$ of its vertical and
horizontal parts.

There is the associated {\em partial connection form}
$\omega$ on $P$ with values in $\fg$, given by
$$ \omega_{p}(\xi)=dL_{p}^{-1}(\xi^{v})\in\fg_{\epsilon(p)}$$
for any $\xi\in \pi^{\ast}(\cF)_{p}$.

\begin{prop}  \label{cpb1}
The partial connection form $\omega$ associated to a connection
$\cH$ on a principal $G$-bundle $P$ has the following properties:
\begin{enumerate}
\item $\omega:\pi^{\ast}(\cF)\ra\fg$ is a map of vector bundles over
      $\epsilon:P\ra G_{0}$,
\item $\omega_{p}\com dL_{p}=\id$,
\item $\Ker(\omega_{p})\subset\Ker((d\epsilon)_{p})$, and
\item $R_{s}^{\ast}\omega=\Ad(\inv\com s)\com \omega$, for any
      local bisection $s:U\ra G_{1}$ of $G$.
\end{enumerate}
\end{prop}
\Rem
Explicitly, the condition (4) means that for any $x\in U$, any
$p\in\epsilon^{-1}(\trg(s(x)))$ and any $\xi\in\pi^{\ast}(\cF)_{p}$
$$ (R_{s}^{\ast}\omega)_{p}(\xi)=\omega_{p s(x)}(dR_{s}(\xi))
   = (dL_{s})_{1_{x}}^{-1}\omega_{p}(\xi) s(x)
   =\Ad(\inv\com s)(\omega_{p}(\xi))\;.$$
Here $dL_{s}$ is the derivative of the diffeomorphism
$L_{s}:\trg^{-1}(U)\ra\trg^{-1}(\trg(s(U)))$
given by $L_{s}(g)=s(\trg(g))g$.
The map $\omega_{p}$ is completely determined by its restriction
to the subspace of those $\xi$ which are in the kernel of $d\epsilon$,
in which case $\anchor(\omega_{p}(\xi))=0$, and the condition (4) may
be expressed simply as
$(R_{s(x)}^{\ast}\omega)(\xi)=\Ad(s(x)^{-1})(\omega(\xi))$.
Conversely, any $\omega$ satisfying the conditions above determines a
connection by
$$ \cH_{p}=\Ker(\omega_{p})\;.$$

\begin{prop}  \label{cpb2}
Let $P$ be a principal $G$-bundle over $B$,
$\cF$ a foliation of $B$ and $\cH$ a flat
$\cF$-partial connection on $P$.
Then each leaf of $\cH$ projects by a covering projection
to a leaf of $\cF$.
\end{prop}
\Proof
Each leaf $\tilde{L}$ of $\cH$ clearly projects by a local
diffeomorphism to a leaf $L$ of $\cF$.
Moreover, $\tilde{L}$ lies in a fiber $\epsilon^{-1}(x)$.
Let $\Iso(\tilde{L})$ be the isotropy group of $\tilde{L}$, i.e.
$$ \Iso(\tilde{L})=\{g\in G(x,x)\,|\,\tilde{L}g=\tilde{L}\}\;.$$
The group $\Iso(\tilde{L})$, equipped with the discrete topology, acts
freely and properly discontinuously on $\tilde{L}$, which shows that
$\pi$ in fact restricts to a covering projection $\tilde{L}\ra L$.
\eop

\section{Lie theory for groupoids}

\subsection{The source-simply connected cover}  \label{sec:scc}

Every finite dimensional Lie algebra is the Lie algebra of a  unique
connected simply connected Lie group. Something similar holds for
{\em integrable} Lie algebroids.

\begin{dfn}  \label{scc1}
A Lie groupoid $G$ is said to be source-connected if $\src^{-1}(x)$
is connected for any $x\in G_{0}$. It is said to be source-simply
connected if each $\src^{-1}(x)$ is connected and simply connected.
\end{dfn}

\begin{ex} \label{scc2} \rm
The monodromy groupoid $\Mon(M,\cF)$ of a foliated manifold
is source-simply connected.
\end{ex}

The following proposition is proved in \cite{M} for the special case
of transitive Lie groupoids.

\begin{prop} \label{scc3}
Let $G$ be a Lie groupoid. There exists a source-simply connected
Lie groupoid $\tilde{G}$ over $G_{0}$
and a morphism of Lie groupoids
$\pi:\tilde{G}\ra G$ over $G_{0}$, inducing an isomorphism
$\cA(\tilde{G})\ra \cA(G)$.
\end{prop}

\Rem
The covering groupoid $\tilde{G}$ is essentially unique, see
Proposition \ref{im1} below.
\vspace{3mm}

\Proof
For each $x\in G_{0}$ let $\src^{-1}(x)^{(0)}$ be the connected component
of $\src^{-1}(x)$ which contains $1_{x}$. It is well-known and easy to see
that the union of these $\src^{-1}(x)^{(0)}$ form a source-connected
open subgroupoid $G^{(0)}$ of $G$ over $G_{0}$ having the same Lie algebroid.
Thus, to prove the proposition we can first replace $G$ by $G^{(0)}$,
and hence assume that $G$ is source-connected.

Now let $\cF$ be the foliation of $G_{1}$ given by the fibers of
$\src$, and let $\Mon(G_{1},\cF)$ be its monodromy groupoid
over $G_{1}$. The space $G_{1}$ is a principal (right) $G$-bundle
over $G_{0}$ (with structure maps $\trg$ for $\pi$ and $\src$ for
$\epsilon$) and this principal action maps leaves to leaves.
Thus $G$ also acts on the monodromy groupoid. By Lemma \ref{ag1}
we can form the quotient Lie groupoid
$\tilde{G}=\Mon(G_{1},\cF)/G$, which is a groupoid over $G_{0}$.
Since any monodromy groupoid is source-simply connected and
$\Mon(G_{1},\cF)$ has the same $\src$-fibers as its quotient
by $G$ by the remark after Lemma \ref{ag1}, $\tilde{G}$ is again
source-simply connected. Finally,
the morphism of Lie groupoids
$$ F:\Mon(G_{1},\cF)\lra G $$
given by $\trg:G_{1}\ra G_{0}$ on objects and by
$F(\sigma)=\sigma(1)\sigma(0)^{-1}$ on arrows
of $\Mon(G_{1},\cF)$ (where $\sigma$ is the homotopy class of a path
inside a leaf of $\cF$) factors to give the required map
$\tilde{G}\ra G$.
\eop

\subsection{Integrability of subalgebroids}  \label{sec:isa}

Let $\fg$ be a Lie algebroid over $M$, and let $N$ be an immersed submanifold
of $M$. A {\em subalgebroid} of $\fg$ over $N$ is a subbundle
$\fh$ of the restriction $\fg|_{N}$ with a Lie algebroid structure such that
the inclusion $\fh\ra \fg$ is a morphism of Lie algebroids.

The same methods involved in the construction of the source-simply connected
groupoid $\tilde{G}$ can be used to prove the following integrability
result. (We point out that the case of transitive groupoids was proved earlier
in \cite{M}, and a local integrability result involving micro-differentiable
groupoids was proved by Almeida, see \cite[p. 158]{M}.)

\begin{prop}  \label{isa1}
Any subalgebroid of an integrable algebroid is integrable.
\end{prop}
\Proof
Consider a Lie groupoid $G$ with associated Lie algebroid $\fg$
and a subalgebroid $\fh$ over $H_{0}$ of $\fg$.
Denote by $I:\fh\ra\fg$ the inclusion
of Lie algebroids over the injective immersion
$\iota:H_{0}\ra G_{0}$, and write
$$ M=H_{0}\times_{G_{0}}G_{1} $$
for the pull-back of $\trg:G_{1}\ra G_{0}$ along $\iota$.
Consider the foliation $\cF$ of $M$
given by
$$ \cF_{(y,g)}=\{(\rho(u),I(u)g)\,|\,u\in\fh_{y}\} $$
The composition of the groupoid $G$ defines on the manifold $M$ the structure of
a principal $G$-bundle over $H_{0}$, and the foliation $\cF$ is preserved by
the $G$-action. Thus the monodromy groupoid $\Mon(M,\cF)$ also
carries a right $G$-action, and by Lemma \ref{ag1} we obtain a quotient groupoid
$$ H=\Mon(M,\cF)/G $$
over $H_{0}$. Its Lie algebroid is easily seen to be isomorphic to $\fh$.
\eop
\vspace{3mm}

\Rem
By the integrability of morphisms between integrable Lie algebroids
(see Proposition \ref{im1} below), there is a map $H\ra G$
which is in fact an immersion. In this sense
the subalgebroid is integrated by an immersed groupoid.

\subsection{Integrability of morphisms between Lie algebroids}  \label{sec:im}

As an application of Proposition \ref{cpb2}, we shall give a quick proof of the
fact that any morphism of integrable Lie algebroids
can be integrated to a unique morphisms of the integral Lie groupoids, provided
that the domain groupoid is source-simply connected. This fact has been proved
earlier by
Mackenzie and Xu \cite{MX}.

\begin{prop}  \label{im1}
Let $G$ and $H$ be Lie groupoids, with $H$ source-simply connected,
and let $\Phi:\fh\ra\fg$ be a morphism of their
Lie algebroids over $\phi:H_{0}\ra G_{0}$.
Then there exists a unique morphism of Lie groupoids
$F:H\ra G$ with $F_{0}=\phi$ which integrates $\Phi$, i.e.
$\cA(F)=\Phi$.
\end{prop}
\Proof
Let $P=H_{1}\times_{G_{0}}G_{1}$ be the pull-back
of $\trg:G_{1}\ra G_{0}$ along the map
$\phi\com\trg:H_{1}\ra G_{0}$. Thus $P$
is a (trivial) principal $G$-bundle over $H_{1}$, with the obvious
right action with respect to the map $\epsilon=\src\com\pr_{2}$.
Let $\cF$ be the foliation of $H_{1}$ by the $\src$-fibers. Now define
a partial connection $\cH$ on $P$ by
$$ \cH_{(h,g)}=\{(v h,\Phi(v) g)\,|\,v\in\fh_{\trg(h)}\}\;.$$
This is indeed a flat connection on $E$ because $\Phi$ preserves the
bracket. Now take any $y\in H_{0}$, and denote by $\tilde{L}_{y}$ the leaf
of $\cH$ through the point $(1_{y},1_{\phi(y)})$.
By Proposition \ref{cpb2}, $\tilde{L}_{y}$ is a covering space over
the corresponding leaf of $\cF$, i.e. the $\src$-fiber $\src^{-1}(y)$.
Since the $\src$-fibers of $H$ are simply connected, the
projection $\tilde{L}_{y}\ra\src^{-1}(y)$
is in fact a diffeomorphism.
Denote by $\nu_{y}$ the inverse of this diffeomorphism.
Now the union of the maps $\nu_{y}$ gives us a map
$\nu:H_{1}\ra P$. Observe that this map is smooth, since it
can be described as an extension of the smooth transversal section
$H_{0}\ra P$ by holonomy. Then we take $F_{1}$ to be the composition
$F_{1}=\pr_{2}\com\nu:H\lra G$. In particular, $F_{1}$ maps
$\src^{-1}(y)$ to $\src^{-1}(\phi(y))$.
It is easy to see that $F_{1}$ together with $F_{0}=\phi$ gives a
morphism of Lie groupoids $F:H\ra G$.
Now for any $v\in\fh_{y}$ we have
$d(F_{1})(v)=d(\pr_{2})(v,\Phi(v))=\Phi(v)$, hence $\cA(F)=\Phi$.
\eop

\subsection{Integrability of actions by Lie groupoids on Lie algebroids} \label{sec:iaga}

Consider a Lie algebroid $\fh$ over $N$. If $q:N\ra B$ is a submersion
which annihilates the anchor of $H$ (i.e. $dq\com\anchor=0$), then
it follows from the Leibniz identity (Section \ref{sec:alg})
that for two sections $Y,Y'\in\Gamma\fh$, the value
$[Y,Y']_{y}$ of the bracket in a point $y\in N$ only depends on the
restrictions of $Y$ and $Y'$ to the fiber $q^{-1}(q(y))$. Thus for
$b\in B$ the fiber $\fh_{b}=\fh|_{q^{-1}(b)}$ is a Lie subalgebroid of $\fh$,
and we can think of $\fh$ as a {\em family} of Lie algebroids over $B$.
For such a family $\fh\ra N\ra B$, the pull-back along any map
$\phi:B'\ra B$ is a Lie algebroid over $B'\times_{B}H_{0}$, which we denote
by $\phi^{\ast}(\fh)$ -- it is a family of Lie algebroids over $B'$.

Now suppose that $G$ is a Lie groupoid, and $\fh$ a family of Lie
algebroids over $G_{0}$, with respect to a surjective submersion
$q:H_{0}\ra G_{0}$ annihilating the anchor of $\fh$, as before.
A left {\em action} of $G$ on $\fh$ {\em along} $q$ is
a morphism of Lie algebroids
$$ \act:\src^{\ast}(\fh)\lra\trg^{\ast}(\fh) $$
over $G_{1}$,
satisfying the unit and cocycle conditions expressing the unit and associativity
laws of an action. Thus, for an arrow $g\in G(x',x)$, the fiber
of $\act$ over $g$ is a Lie algebroid morphism
$$ \act_{g}:\fh_{x'}\lra \fh_{x} \;\;\;\;\;\;\;\;\;\; \act_{g}(u)=gu \;$$
and the unit and cocycle conditions state that $\act_{1_{x}}=\id$ and
$\act_{g}\act_{g'}=\act_{gg'}$.
In particular, the effect of the map $\act$
on the base manifolds of the algebroids
is a map $G_{1}\times_{G_{0}}H_{0}\ra H_{0}\times_{G_{0}}G_{1}$ over $G_{1}$,
and its composition with $\pr_{1}$ defines a $G$-action on the space
$H_{0}$. For $g\in G(x',x)$ and $y\in q^{-1}(x')$ this action is denoted
by $gy\in q^{-1}(x)$ as before. Thus if $u\in\fh_{y}$ is a point in the
algebroid $\fh$ over $y$, then $gu\in\fh_{gy}$.

For example, a left action by a Lie groupoid $G$ on a Lie
groupoid $H$ induces an action of $G$ on the Lie algebroid
$\fh$ of $H$, by taking the
derivative of the action on arrows.

We will now show that if $\fh$ is integrable, then so is any action
by $G$ on $\fh$, i.e. any such action comes from an action by $G$ on
the source-simply connected integral $H$ of $\fh$.

\begin{theo}  \label{iaga1}
Let $\fh$ be the Lie algebroid of an source-simply connected
Lie groupoid $H$. Any action of a Lie groupoid $G$ on the Lie algebroid $\fh$
can be integrated uniquely to an action of $G$ on $H$.
\end{theo}
\Proof
As before, we denote by $q:H_{0}\ra G_{0}$ the structure map of this action.
Note that since $dq\com\anchor=0$, the map $q\com\trg$ is locally constant
on the $\src$-fibers of $H$, and hence it is constant since $H$ is source-connected.
Thus $q\com\src=q\com\trg$. Now consider the pull-back manifold
$$ M=G_{1}\times_{G_{0}}H_{1}=\{(g,h)\,|\,\src(g)=q(\src(h))\}\;.$$
This manifold is equipped with a foliation $\cF$, whose leaves are the fibers of
the map $\id\times\src:M=G_{1}\times_{G_{0}}H_{1}\ra G_{1}\times_{G_{0}}H_{0}$.
Over $M$ there is a principal (right) $H$-bundle
$$ P=G_{1}\times_{G_{0}}H_{1}\times_{H_{0}}H_{1}
    =\{(g,h,h')\,|\,(g,h)\in M,\, g\trg(h)=\trg(h') \} $$
with structure maps $\src\com\pr_{3}:P\ra H_{0}$ and action
$(g,h,h')h''=(g,h,h'h'')$. The foliation $\cF$ of $M$ lifts
along the projection $\pi:P\cong M\times_{H_{0}}H_{1}\ra M$
to a foliation $\cH$
of $P$, whose tangent space $\cH_{(g,h,h')}\subset T_{(g,h,h')}(P)$ is defined
in terms of the action of $G$ on $\fh$ by
$$ \cH_{(g,h,h')}=\{ (0g,uh,(gu)h') \,|\, u\in\fh_{\trg(h)} \} \;.$$
To explain this notation, let $g\in G(x',x)$ and
$y=\trg(h)$, $gy=\trg(h')$. Then $0g$ is the zero tangent vector in
$T_{g}(G_{1})$, and $uh\in T^{\src}_{h}(H_{1})$.
Furthermore, $gu\in \fh_{gy}$ and hence
$(gu)h'\in T^{\src}_{h'}(H_{1})$.

Note that, since the action of $G$ on $\fh$ preserves the bracket, this subbundle
$\cH$ of $T(P)$ is integrable, so that $\cH$ is indeed a foliation of $P$.
It is clear from this definition that $\cH$ is a partial flat connection
on the principal $H$-bundle $P$ over the foliation $(M,\cF)$.
So by Proposition \ref{cpb2} and the fact that the leaves of $\cF$ are simply
connected, the projection $\pi:P\ra M$ of the principal bundle restricts
to a diffeomorphism $L\ra \pi(L)$ from any leaf $L$ of $\cH$ to the leaf
$\pi(L)$ of $\cF$. Now consider the complete transversal sections
$S$ of $(M,\cF)$ and $T$ of $(P,\cH)$, defined by
$$ S=\{ (g,1_{y}) \,|\, g\in G_{1},\,y\in H_{0},\,\src(g)=q(y)\} $$
and
$$ T=\{ (g,1_{y},1_{y'}) \,|\, g\in G_{1},\,y,y'\in H_{0},\,
     \src(g)=q(y),\,gy=y'\} \;.$$
Let $\phi:M\ra P$ be the unique section of $\pi$ which sends $T$ into $S$
by $\phi(g,1_{y})=(g,1_{y},1_{q(y)})$ and which maps the leaf through
$(g,1_{y})$ to the leaf through $\phi(g,1_{y})$. Then
$$ \act=\pr_{3}\com\phi:G_{1}\times_{G_{0}}H_{1}\lra H_{1} $$
is the action map which integrates the given action of $G$ on $\fh$.
Indeed, to see that it respects the composition in $G$ and in $H$, note that
for a fixed $g\in G(x',x)$ the map $\act_{g}=\act(g,\oo)$ is the unique
Lie groupoid morphism $H_{q^{-1}(x')}\ra H_{q^{-1}(x)}$, integrating the
algebroid map $\fh_{x'}\ra \fh_{x}$ which sends $u$ to $gu$
(see Proposition \ref{im1}).
\eop

\section{Derivations and infinitesimal actions}

\subsection{Derivations on Lie algebroids} \label{sec:dla}

Let $M$ be a manifold , and $\fX(M)$ the Lie algebra of vector fields on $M$.
Recall that a derivation on $\fX(M)$ is an $\RR$-linear map
$D:\fX(M)\ra\fX(M)$ satisfying
$$ D([X,X'])=[D(X),X']+[X,D(X')] $$
for any two vector fields $X,X'\in\fX(M)$.
Each vector field $V$ on $M$ defines an (inner) derivation
$D$ on $\fX(M)$ given by $D(X)=[V,X]$, and in fact any
derivation on $\fX(M)$ is of this form for a unique vector field
$V\in\fX(M)$ \cite{T}.
In this section we will prove a similar result for derivations on the Lie
algebra of sections of a Lie algebroid.

\begin{dfn}  \label{dla1}
A derivation on a Lie algebroid $\fg$ over $M$ is a pair
$(D,V)$ consisting of an $\RR$-linear map $D:\Gamma\fg\ra\Gamma\fg$ and
a vector field $V$ on $M$ such that
\begin{enumerate}
\item [(i)]    $D([X,X'])=[D(X), X']+[X,D(X')]$,
\item [(ii)]   $D(fX)=fD(X)+V(f)X$ and
\item [(iii)]  $\anchor(D(X))=[V,\anchor(X)]$
\end{enumerate}
for any $X,X'\in\Gamma\fg$ and $f\in\eC(M)$.
\end{dfn}
\Rem
If the rank of
$\fg$ is not zero, the properties (i) and (ii) imply
property (iii).
The vector space of all derivations on $\fg$, denoted by
$\Der(\fg)$, is a Lie algebra with respect to the bracket
$[(D,V),(D',V')]=(D\com D'-D'\com D,[V,V'])$.

\begin{exs} \rm \label{dla2}
(i) For a Lie algebra (viewed as a Lie algebroid over a one point space)
we recover the usual notion of a derivation.

(ii) Suppose that $E\ra M$ is a vector bundle, viewed as a Lie algebroid
as in Example \ref{alg1} (iii). A derivation on this Lie algebroid
consists of a vector field $V$ on $M$ and a partial connection
(a covariant differential operator) $D=\nabla_{V}$ on $E$.

(iii) Let $\cF$ be a foliation of a manifold $M$,
viewed as a Lie algebroid over $M$ with injective anchor map.
The Lie algebra of derivations $\Der(\cF)$ can be identified
with the Lie algebra $L(M,\cF)$ of projectable vector fields
on $(M,\cF)$ (see \cite{Mo}).
\end{exs}

For a Lie groupoid $G$, there is an associated tangent Lie groupoid
$T(G)$ over $T(G_{0})$. Its manifold of arrows is $T(G_{1})$, while the source
and the target maps $T(G_{1})\ra T(G_{0})$ and the multiplication map
$T(G_{1})\times_{T(G_{0})}T(G_{1})\cong T(G_{1}\times_{G_{0}} G_{1})\ra T(G_{1})$
are the derivatives of those of $G$. The bundle projections
$\pi_{i}:T(G_{i})\ra G_{i}$, $i=0,1$, define a morphism of Lie groupoids
$$ T(G)\lra G \;.$$

\begin{dfn}  \label{dla3}
A multiplicative vector field on a Lie groupoid $G$
is a morphism of Lie groupoids $W:G\ra T(G)$, which is a section of
the projection $T(G)\ra G$.
\end{dfn}
\Rem
Multiplicative vector fields were studied in \cite{MX1}.
A multiplicative vector field $W$ on $G$ is, in other words, a pair
of vector fields $W_{0}$ on $G_{0}$ and $W_{1}$ on $G_{1}$, such that
$W_{1}$ is projectable to $W_{0}$ along both $\src$ and $\trg$,
$W_{0}$ is projectable to $W_{1}$ along $\uni:G_{0}\ra G_{1}$,
and $(W_{1},W_{1})\in \fX(G_{1}\times_{G_{0}}G_{1})$ is projectable to $W_{1}$
along the multiplication of $G$.
The last two conditions mean that
$$ (W_{1})_{1_{x}}=1_{(W_{0})_{x}}=d(\uni)((W_{0})_{x}) \;\;\;\;\;\;\;\;
   \mathrm{and} \;\;\;\;\;\;\;\;
   (W_{1})_{gg'}=(W_{1})_{g}(W_{1})_{g'}\;, $$
the latter composition being the one in $T(G)$.
Note that $W_{0}$ is determined by $W_{1}$, in fact the restriction of
$W_{1}$ to $G_{0}$ is tangent to $G_{0}$ and can be identified with
$W_{0}$. Therefore we will simply write $W$ for $W_{1}$ and identify $W_{0}$ with
$W_{1}|_{G_{0}}$.

The Lie brackets of vector fields on $G_{0}$ and $G_{1}$ respect projectability
along the maps between $G_{0}$ and $G_{1}$ and hence define a Lie bracket on
multiplicative vector fields on $G$. In this way, the multiplicative vector fields on
$G$ form a Lie algebra, denoted by $\fX^{\mlt}(G)$.

\begin{lem}  \label{dla4}
For any Lie groupoid $G$ we have
$[\fX^{\mlt}(G),\fX^{\src}_{inv}(G)]\subset \fX^{\src}_{inv}(G)$.
\end{lem}
\Proof
Let $W\in\fX^{\mlt}(G)$ and $X\in\fX^{\src}_{inv}(G)$.
The invariance of $X$ is equivalent to the condition that the vector
field $(X,0)$ on $G_{1}\times_{G_{0}}G_{1}$ is projectable to $X$
along $\mlt$. Since $W$ is projectable along $\src$ and
$X$ is tangent to the fibers of $\src$, $[W,X]$ is tangent to the fibers
of $\src$ as well. Since both $(W,W)$ and $(X,0)$ are projectable along $\mlt$,
so is $[(W,W),(X,0)]=([W,X],0)$, and
$d\mlt\com ([W,X],0)=[W,X]\com\mlt$.
\eop
\vspace{3mm}

Let $G$ be a Lie groupoid and $\fg$ the Lie algebroid
associated to $G$. Any multiplicative vector field $W$ on $G$
gives us a derivation $\cL(W)=(\cL_{W},W|_{G_{0}})$
on $\fg$ by
$$ \cL_{W}(X)=[W,X] \;\;\;\;\;\;\;\;\;\; X\in\Gamma\fg\cong\fX^{\src}_{inv}(G)\;.$$
Indeed, Lemma \ref{dla3} implies that the image of $\cL_{W}$ is
in $\Gamma\fg$, while it is easy to check that the properties (i), (ii) and (iii)
in Definition \ref{dla1} are satisfied. Moreover, $\cL$ is a morphism of Lie algebras
$$ \cL:\fX^{\mlt}(G)\lra\Der(\fg)\;.$$
The following theorem follows from the results of \cite{MX1}:

\begin{theo}  \label{dla5}
If $G$ is a source-simply connected Lie groupoid, then the map
$\cL$ is an isomorphism of Lie algebras.
\end{theo}

\Proof
First we will show that $\cL$ is injective. Let $W\in\Ker\cL$. In particular,
$d\src\com W=d\trg\com W=0$, $W|_{G_{0}}=0$ and $[W,\fX^{\src}_{inv}(G)]=0$.
Take any $g\in G_{1}$, and consider the fiber $G(\src(g),\oo)$.
Since any tangent vector on $G(\src(g),\oo)$ can be extended to an
$\src$-vertical invariant vector field, the condition
$[W,\fX^{\src}_{inv}(G)]=0$ implies that the subset of zeros of
$W$ is open in $G(\src(g),\oo)$. But $G(\src(g),\oo)$ is connected
and $W_{\src(g)}=0$, thus $W_{g}=0$.

Next we will prove that $\cL$ is surjective. Let $(D,V)$ be
a derivation on the Lie algebroid $\fg$ associated to $G$.
Take any $u\in\fg_{x}$. For any $v\in\fg_{x}$ let
$\tau(u,v)\in\Ker((d\pi)_{u})\subset T_{u}(\fg)$ be given
by
$$ \tau(u,v)(f)=\left.\frac{d}{dt}\right|_{t=0}\!\!\!\!\! f(u+tv)
   \;\;\;\;\;\;\;\;\;\; f\in\eC(\fg)\;.$$
In other words, the map
$\tau(u,\oo)$ is the natural isomorphism between
$\fg_{x}$ and $\Ker((d\pi)_{u})$. Now define
$$ \Xi(u)=(dX)_{x}(V_{x})-\tau(u,D(X)_{x})
   \in T_{u}(\fg)\;,$$
for any section $X\in\Gamma\fg$ satisfying
$X_{x}=u$. The property (ii) of $D$
(Definition \ref{dla1}) implies that the definition of $\Xi(u)$
does not depend on the choice of $X$, so we get a map
$$ \Xi:\fg\lra T(\fg)\;.$$
Moreover, $\Xi$ is a bundle map over $V$
by \cite[Proposition 2.5]{MX1}. Now recall from \cite{MX2}
that $T(\fg)$ has a natural structure of a Lie algebroid over $T(G_{0})$
and that there is a natural isomorphism of Lie algebroids
$j:T(\fg)\ra\cA(T(G))$ over $T(G_{0})$.
It follows from \cite[Theorem 4.4]{MX1}
that $\Xi$ is a morphism of Lie algebroids, hence
$\tilde{\Xi}=j\com\Xi:\fg\ra\cA(T(G))$ is a morphism of
Lie algebroids over $V$ as well. Since $G$ is source-simply connected,
Proposition \ref{im1} implies that
$\tilde{\Xi}$ can be integrated to a unique morphism of Lie
groupoids
$$ W:G\lra T(G)\;.$$
Note that $W$ is a multiplicative vector field on $G$
because $\tilde{\Xi}$ is a
section of the projection $\cA(T(G))\ra\fg$, and that $W|_{G_{0}}=V$.
Finally, \cite[Theorem 3.9]{MX1} implies that $\cL_{W}=D$.
\eop

\subsection{Infinitesimal actions} \label{sec:ia}

Let $\fh$ be a Lie algebroid over $N$ and $q:N\ra M$
a surjective submersion. Thus the Lie algebra $\Gamma\fh$ has
the structure of a $\eC(M)$-module induced by the composition
with $q$. Assume that $dq\com\anchor=0$, i.e. that
$\fh$ is a family of Lie algebroids over $M$. This implies
that $\anchor(Y)(f\com q)=0$ for any $Y\in\Gamma\fh$ and any
$f\in\eC(M)$, and hence the Lie bracket on $\Gamma\fh$ is
$\eC(M)$-bilinear. Furthermore,
the Lie algebra of derivations $\Der(\fh)$ on $\fh$ becomes
a $\eC(M)$-module via $f(D,V)=((f\com q)D,(f\com q)V)$.
If $H$ is a source-simply connected integral of $\fh$, the Lie
algebra $\fX^{\mlt}(H)$ of multiplicative
vector fields on $H$ is also
a $\eC(M)$-module and the isomorphism $\cL$ of Theorem \ref{dla5}
is $\eC(M)$-linear.

\begin{dfn} \label{ia1}
Let $\fg$ and $\fh$ be Lie algebroids over $M$ respectively $N$,
and let $q:N\ra M$ be
a surjective submersion such that $dq\com\anchor=0$.
An (infinitesimal) action of $\fg$ on $\fh$ along $q$
is a homomorphism of Lie algebras $\nabla:\Gamma\fg\ra\Der(\fh)$,
$\nabla(X)=(\nabla_{X},R(X))$, which is $\eC(M)$-linear
and for which each $R(X)$ is projectable to $\anchor(X)$ along $q$.
\end{dfn}
\Rem
In particular, $R:\Gamma\fg\ra\fX(N)$ is a homomorphism of Lie
algebras satisfying $\nabla_{X}(f'Y)=f'\nabla_{X}Y+R(X)(f')Y$
for any $X\in\Gamma\fg$, $Y\in\Gamma\fh$ and
$f'\in\eC(N)$. Furthermore, the $\eC(M)$-linearity of $\nabla$
implies that $\nabla_{fX}(Y)=(f\com q)\nabla_{X}(Y)$ and
$R(fX)=(f\com q)R(X)$.
The second equality in fact follows from the first if
the rank of $\fh$ is not zero.
The projectability of $R(X)$ to $\anchor(X)$ along $q$
means that $R(X)(f\com q)=\anchor(X)(f)\com q$
for any $f\in\eC(M)$ and $X\in\Gamma\fg$.
We shall say that
$\nabla$ is an action {\em over} $R$.
Our definition is clearly equivalent with
\cite[Definition 3.6]{HM}.

\begin{exs} \rm \label{ia2}
(i) If $\fg$ and $\fh$ are Lie algebras, we recover the usual notion
of an action \cite{S}.

(ii) If $\fg$ is a Lie algebra and $\fh$ is the tangent bundle of a manifold
$N$, an infinitesimal action of $\fg$ on $\fh$ is the same thing as an infinitesimal
action of $\fg$ on $N$ (Example \ref{alg1} (vi)).
More generally, if $\fh$ is a foliation $\cF$ of a manifold $N$
(Example \ref{alg1} (v)),
an infinitesimal action of $\fg$ on $\cF$ is a Lie algebra map
$\fg\ra L(N,\cF)$ into the projectable vector fields on $(N,\cF)$.

(iii) Let $E\ra M$ be a vector bundle over a foliated manifold
$(M,\cF)$. We can consider $\cF$ and $E$ as Lie algebroids
over the same manifold $M$ (Examples \ref{alg1} (iv) and (v)),
and an action of $\cF$ on $E$
along the identity map is the same thing as an affine flat $\cF$-partial
connection on $E$.
\end{exs}

Let $\nabla$ be an action of $\fg$ on $\fh$ along $q:N\ra M$
over $R$.
Then the {\em semi-direct product} \cite{HM}
of $\fg$ and $\fh$ with
respect to $\nabla$ is a Lie algebroid $\fg\ltimes\fh$ over $N$
given as follows: as a vector bundle it is the direct sum
$q^{\ast}\fg\oplus\fh$, the anchor is given by
$$ \anchor(q^{\ast}X\oplus Y)=R(X)+\anchor(Y)\;,$$
and the bracket by
$$ [q^{\ast}X\oplus Y,q^{\ast}X'\oplus Y']
   =q^{\ast}[X,X']\oplus([Y,Y']+\nabla_{X}(Y')-\nabla_{X'}(Y) )\;. $$
Here $q^{\ast}X=(\id,X\com q)\in\Gamma q^{\ast}\fg$.
Since the sections of this form span $\Gamma q^{\ast}\fg$ as a
$\eC(N)$-module,
we can extend the definition of the anchor and
of the bracket to all the sections of $q^{\ast}\fg\oplus\fh$
by the $\eC(N)$-linearity of the anchor and by the Leibniz identity.

Observe that we have an exact sequence
$$ 0 \lra \fh \stackrel{j}{\lra} \fg\ltimes\fh
   \stackrel{\pi}{\lra} q^{\ast}\fg \lra 0 $$
of vector bundles over $N$. The action of $\fg$ on $\fh$ can be recovered
from this sequence by
$$ j(\nabla_{X}(Y))=[(q^{\ast}X,0),(0,Y)]\;.$$
Semi-direct products are related to split exact sequences in the usual way.
Explicitly, suppose that $\fk$ is a Lie algebroid over $N$
which fits into an exact sequence
$$ 0 \lra \fh \stackrel{j}{\lra} \fk
   \stackrel{\pi}{\lra} q^{\ast}\fg \lra 0 $$
of vector bundles over $N$. Suppose that $j$ is a map of Lie
algebroids over $N$, and that $\pi$ is given by map of Lie algebroids
$\fk\ra\fg$ over $q$.
Consider a splitting of this exact sequence by a map
$i:q^{\ast}\fg\ra\fk$ of vector bundles. The map $i$ defines a
``connection'' $\nabla:\Gamma\fg\otimes\Gamma\fh\ra\Gamma\fh$,
$\nabla(X\otimes Y)=\nabla_{X}(Y)$,
by
$$ j(\nabla_{X}(Y))=[i(q^{\ast}X),j(Y)]_{(\fk)}\;.$$
The curvature $2$-form of this connection is the map
$\kappa:\Gamma\fg\wedge\Gamma\fg\ra\Gamma\fh$ given by
$$ j(\kappa(X,X'))=[i(q^{\ast}X),i(q^{\ast}X')]_{(\fk)}
   -i(q^{\ast}[X,X']_{(\fg)})\;.$$
The connection $\nabla$ is {\em flat} if
$[\kappa(X,X'),Y]=0$ for every $Y\in\Gamma\fh$.
In particular, this
is the case if $\kappa=0$, i.e. if $i$ preserves the bracket.
In this case, $\nabla$ is an action
of $\fg$ on $\fh$ along $q$ as defined above,
with $R(X)=\anchor(i(q^{\ast}X))$, and
$\fk$ is isomorphic to $\fg\ltimes\fh$.

\begin{exs} \rm \label{ia3}
(i) Let $P$ be a principal $H$-bundle over $M$ for a Lie group $H$, and let
$G$ be the gauge groupoid over $M$ (Example \ref{gpd1} (vii)),
with Lie algebroid $\fg$. There is an exact ``Atiyah'' sequence over $M$,
$$ 0\lra\fh^{tw} \lra \fg \lra T(M) \lra 0 \;,$$
where $\fh^{tw}$ is the bundle of Lie algebras over $M$ obtained by
twisting the trivial bundle $M\times\fh$ by the adjoint action
of the cocycle defining the principal bundle \cite{A}. Here $\fh$ is the
Lie algebra associated to $H$.
An Ehresmann connection is the same thing as a splitting of this exact sequence.
Flat connections correspond to actions of $T(M)$ on $\fh^{tw}$,
and represent $\fg$
as a semi-direct product $T(M)\ltimes \fh^{tw}$.

(ii)  If $\fg$ and $\fh$ are Lie algebroids with injective anchor maps,
then the semi-direct product $\fg\ltimes\fh$ of an infinitesimal action
of $\fg$ on $\fh$ along $q:N\ra M$
again has injective anchor. (Indeed, if $X\in\Gamma\fg$ and
$Y\in\Gamma\fh$ are such that $(R(X)+\anchor(Y))_{y}=0$ for a point $y\in N$,
then $0=dq(R(X)+\anchor(Y))_{y}=\anchor(X)_{q(y)}$, whence $X_{q(y)}=0$. Then also
$R(X)_{y}=0$ because $R$ is $\eC(M)$-linear, whence $\anchor(Y)_{y}=0$ so
$Y_{y}=0$.) Thus, a semi-direct product of two foliations is again a foliation.
\end{exs}

\section{Integrability of semi-direct products}

\subsection{Infinitesimal actions on foliations}  \label{sec:iaf}

In this section we consider infinitesimal actions of a Lie algebroid $\fg$ on
a Lie algebroid $\fh$ in the special case where $\fh$ is a foliation $\cF$ of
a manifold $N$ (Example \ref{alg1} (v)). As noted before, such an algebroid
$\cF$ is always integrable, e.g. by the monodromy groupoid of the foliation.
We will prove the following result:

\begin{theo}  \label{iaf1}
For any action of an integrable Lie algebroid $\fg$ on a foliation $\cF$, the
semi-direct product $\fg\ltimes\cF$ is  integrable.
\end{theo}

\Proof
Let $G$ be a Lie groupoid which integrates $\fg$.
First, let us spell out what it means for the Lie algebroid
$\fg$ of $G$ to act on a foliation $\cF$ of a manifold $N$.
First, we have a surjective submersion $q:N\ra G_{0}$, and the leaves
of $\cF$ are contained in the fibers of $q$. Thus, any fiber
$q^{-1}(x)$ is itself a foliated manifold. Next, there is a $\eC(G_{0})$-linear
Lie algebra map $R:\Gamma\fg\ra L(N,\cF)$ into the projectable vector fields
on $(N,\cF)$. In particular, any $X\in\Gamma\fg$ induces a derivation
$\nabla_{X}$ on the vector fields $Y$ on $N$ which are tangent to $\cF$ by
$\nabla_{X}(Y)=[R(X),Y]$.
Finally, this map satisfies the condition that $R(X)$ is projectable along
$q$ to the anchor $\anchor(X)$, i.e.
$dq((R(X)_{y})=\anchor(X)_{q(y)}$ for any $y\in N$.

Now consider  the principal $G$-bundle $P=N\times_{G_{0}}G_{1}$ over $N$.
Here $P$ is the pull-back of $\trg:G_{1}\ra G_{0}$ along $q:N\ra G_{0}$, with
the evident right $G$-action $(y,g)g'=(y,gg')$.
Consider on $P$ the foliation $\cG$ whose tangent space
$\cG_{(y,g)}\subset T_{(y,g)}(P)$ consists of pairs
$$ (R(X)_{y}+Y_{y},X_{\trg(g)}g)\;,$$
where $X\in\Gamma\fg$ and $Y$ is a vector field on $N$ tangent to the
foliation $\cF$.
In other words, the sections of $\cG$ are spanned as a $\eC(P)$-module by
sections of the form $(R(X)+Y,X)$. Here we denote the invariant vector field
on $G_{1}$ associated to $X$ again by $X$.
The subbundle $\cG$ is involutive;
indeed, using the fact that $R$ preserves the bracket, we have
\begin{eqnarray}
[(R(X)+Y,X),(R(X')+Y',X')]
= (R([X,X'])+Z,[X,X']) \label{eq1} \;,
\end{eqnarray}
where $Z=\nabla_{X}(Y')-\nabla_{X'}(Y)+[Y,Y']$ is again tangent to the
foliation $\cF$.

The leaves of this foliation $\cG$ are contained in the fibers of the map
$\src\com\pr_{2}:P\ra G_{0}$ which is a part of the principal bundle structure,
and the $G$-action maps leaves into leaves. Thus $G$ also acts on the Lie
groupoid $\Mon(P,\cG)$. By Lemma \ref{ag1}, the quotient $K=\Mon(P,\cG)/G$ is
a Lie groupoid over $N$. We claim that $K$ integrates the semi-direct product
$\fg\ltimes\cF$. Indeed, the map
$$ \Phi:\cG\lra\fg\ltimes\cF\;,$$
sending a $\cG$-tangent vector field
$(R(X)+Y,X)$ to $q^{\ast}X\oplus Y$,
preserves the bracket (compare with the Equation (\ref{eq1})).
Note also that $\Phi$ respects the anchor because
$$ d(\pr_{1})(R(X)_{y}+Y_{y},X_{\trg(g)}g)
   =R(X)_{y}+Y_{y}=\anchor(q^{\ast}X\oplus Y)_{y}\;,$$
so it is a morphism of Lie algebroids over $\phi=\pr_{1}:P\ra N$.
This map $\Phi$ induces an isomorphism of the fibers
$\cG_{(y,g)}\cong(\fg\ltimes\cF)_{y}$.
Since the $\src$-fibers of the quotient $K$ are the same
as those of $\Mon(P,\cG)$, the vector space
$\cG_{(y,g)}$ is isomorphic to the fiber $\fk_{y}$ of the Lie
algebroid $\fk$ of $K$, hence
$\Phi$ induces an isomorphism
$\fk\ra\fg\ltimes\cF$.
\eop

\begin{ex} \rm \label{iaf2}
If $\fg$ is the Lie algebra of a Lie group, the algebroid
associated to an infinitesimal action of $\fg$ on a manifold
(Example \ref{alg1} (vi)) is always integrable.
This result is due to Dazord \cite{D} (see also \cite{Pa}),
and is a special case of the previous theorem.
\end{ex}

\subsection{Infinitesimal actions along proper maps} \label{sec:iapm}

Consider  an infinitesimal action of a Lie algebroid $\fg$
over $G_{0}$ on a manifold $N$ (viewed as the trivial Lie algebroid over $N$).
Theorem \ref{iaf1} implies that the semi-direct product $\fg\ltimes N$
is integrable by some Lie groupoid
whenever $\fg$ is integrable. However, if $G$ integrates $\fg$,
one can in general not integrate the infinitesimal action
to an action of $G$ on $N$. In this section, we will first show that
such an action can be integrated in the special case where the map $q:N\ra G_{0}$
(which is a part of the infinitesimal action) is proper.

\begin{theo} \label{iapm1}
Let $G$ be a source-simply connected Lie groupoid and suppose that
the Lie algebroid $\fg$ of $G$
acts on a manifold $N$ along a proper map $q:N\ra G_{0}$.
Then there exists an action of $G$ on $N$ along $q$ which integrates
the infinitesimal action in the sense that
$\fg\ltimes N$ is isomorphic to the Lie algebroid of the translation
groupoid $G\ltimes N$.
\end{theo}

\Proof
We view $N$ as a manifold foliated by points, and consider (as in the proof
of Theorem \ref{iaf1}) the foliation $\cG$ on $N\times_{G_{0}}G_{1}$, whose tangent
space $\cG_{(y,g)}\subset T_{(y,g)}(N\times_{G_{0}}G_{1})$ consists of pairs
$((R(X)_{y},X_{\trg(g)}g)$ for $X\in\Gamma\fg$.
Also, we consider the foliation $\cF$ of $G_{1}$ by the $\src$-fibers.

It is clear that the projection $\pi=\pr_{2}:N\times_{G_{0}}G_{1}\ra G_{1}$ maps
$\cG$ to $\cF$ and restricts to a local diffeomorphism from any leaf $L$ of
$\cG$ to a leaf $L'$ of $\cF$.
The projection $\pi$ is also proper, because it is a pull-back of $q$ which is proper.
We will show that this implies that $\pi(L)=L'$ and that $\pi|_{L}:L\ra L'$ is a covering
projection. To see this, take any arrow $g\in L'$ and choose vector fields
$X_{1},\ldots,X_{k}$ on $G_{1}$ such that their values at $g$ form a basis of
$\cF_{g}$ (in fact, we can choose $X_{1},\ldots,X_{k}$ to be in $\fX^{\src}_{inv}(G)$).
Suppose that $U$ is a small open neighbourhood of $g$ in $L'$ and $\varepsilon>0$
such that the local flow $\varphi^{i}_{t}:U\ra L'$ of $X_{i}$ is well-defined
for any $t\in (-\varepsilon,\varepsilon)$, $i=1,\ldots,k$. We can
choose $\varepsilon$ so small that the local flows give us an open embedding
$\psi:(-\varepsilon,\varepsilon)^{k}\ra L'$ by
$$ \psi(t_{1},\ldots,t_{k})=(\varphi^{1}_{t_{1}}\com\ldots\com\varphi^{k}_{t_{k}})(g)\;.$$
Let $\tilde{X}_{i}$ be the unique vector field on $N\times_{G_{0}}G_{1}$ tangent to
$\cG$ which projects to $X_{i}$ along $\pi$.
Since $q$ is proper we can take $U$ and $\epsilon$ so small
that the local flow
$\tilde{\varphi}^{i}_{t}:\pi^{-1}(U)\ra \pi^{-1}(L')$ of $\tilde{X}_{i}$ is well-defined
for any $t\in (-\epsilon,\epsilon)$, $i=1,\ldots,k$.
Note that $\pi\com\tilde{\varphi}^{i}_{t}=\varphi^{i}_{t}\com\pi$ because
$\tilde{X}_{i}$ is a lift of $X_{i}$. In particular, the map
$\tilde{\psi}:(-\varepsilon,\varepsilon)^{k}\times\pi^{-1}(g)\ra \pi^{-1}(L')$ given by
$$ \tilde{\psi}(t_{1},\ldots,t_{k},(y,g))=
   (\tilde{\varphi}^{1}_{t_{1}}\com\ldots\com\tilde{\varphi}^{k}_{t_{k}})(y,g) $$
is well-defined and satisfies $\pi\com\tilde{\psi}=\psi\com\pr_{1}$.
It clearly follows that $\tilde{\psi}$ is an open embedding.
The map $\tilde{\psi}$ also maps the
product foliation of $(-\varepsilon,\varepsilon)^{k}\times\pi^{-1}(g)$
to $\cG$.
In particular, the projection $\pi|_{L}:L\ra L'$ is a covering.

Now recall that in the case at hand
$L'$ is an $\src$-fiber of $G$ and hence simply connected.
It follows that the covering projection $\pi|_{L}:L\ra L'$ is in fact
a diffeomorphism.
We can now define a $G$-action on $N$
$$ \act:G_{1}\times_{G_{0}}N\lra N $$
along $q$ as follows: for any arrow $g$ of $G$ and any point $y$ of $N$
satisfying $\src(g)=q(y)$, let
$gy$ be the unique point of $N$ such that
$(gy,g)$ lies on the same leaf of $\cG$ as $(y,1_{q(y)})$.
It is clear that the map $\act$
satisfies the identities for an action.
To see that $\act$ is smooth, observe that a lift of a holonomy extension of a
path $\gamma$ in a leaf of $\cF$ can be obtained as a holonomy extension
of a lift of $\gamma$.
It is straightforward to check that this
$G$-action on $N$ indeed integrates
the infinitesimal action of $\fg$ on $N$.
\eop
\vspace{3mm}

\Rem
Note that the assumption that $q$ is proper can be replaced by the assumption
that all the vector fields $R(X)$ are complete.

\begin{cor}  \label{iapm2}
Let $G$ and $H$ be source-simply connected Lie groupoids with Lie algebroids $\fg$
respectively $\fh$, and suppose that
$\fg$ acts on $\fh$ along a proper map $q:H_{0}\ra G_{0}$. Then
there exists an action of $G$ on $H$ which integrates the infinitesimal action,
in the sense that the semi-direct product $\fg\ltimes\fh$ is isomorphic to the
Lie algebroid of the semi-direct product groupoid $G\ltimes H$.
\end{cor}

\Proof
Consider the action of $G$ on $H_{0}$ given by Theorem \ref{iapm1}.
The Lie algebroid $\fh$ is a family of Lie algebroids over $G_{0}$,
so its pull-back along $\pr_{2}:G_{1}\times_{G_{0}}H_{0}\ra H_{0}$
is again a Lie algebroid. This pull-back is equipped with a foliation whose leaves are
the fibers of $\pr_{2}$. These fibers are isomorphic to the $\src$-fibers of $G$.
The action of $\fg$ on $\fh$ defines a flat partial connection
on $\pr_{2}^{\ast}\fh$ along the leaves of the foliation of $G_{1}\times_{G_{0}}H_{0}$.
Since the leaves are simply connected, there is a well-defined transport
along the leaves, which defines an action of $G$ on the algebroid $\fh$.
This action can now be integrated by Theorem \ref{iaga1}.
Further details are straightforward.
\eop
\vspace{3mm}

Let $\fg$ be a Lie algebroid over $M$, and let $\fh$ be a bundle
of Lie algebras over $M$ (Example \ref{alg1} (iv)). An {\em extension} $\fk$
of $\fg$ by $\fh$ is an exact sequence of Lie algebroids over $M$
of the form
$$ 0 \lra \fh \lra \fk \lra \fg \lra 0\;.$$
It is called {\em split} if it splits by a morphism of Lie algebroids
$\fg\ra\fk$ over $M$.

\begin{cor}  \label{iapm3}
Any split extension of an integrable Lie algebroid over $M$
is integrable.
\end{cor}

\Proof
A split extension as above represents $\fk$ as the semi-direct product
$\fg\ltimes\fh$ for an action by $\fg$ on $\fh$ along the
identity $M\ra M$. Thus the result follows from the previous corollary and
integrability of $\fh$ (Example \ref{lalg2} (iii)).
\eop

\subsection{Infinitesimal actions on the Lie algebroids of
source-compact groupoids}  \label{sec:iascg}

We consider again two Lie groupoids $G$ and $H$
with Lie algebroids $\fg$ and $\fh$,
respectively, and an infinitesimal action of $\fg$ on $\fh$ along a map
$q:H_{0}\ra G_{0}$. In this section we shall prove that the semi-direct product
$\fg\ltimes\fh$ is integrable in the case where $H$ is a source-compact
source-simply connected groupoid.

\begin{dfn}  \label{iascg0}
A Lie groupoid $H$ is called source-compact if it is Hausdorff
and if the source map $\src:H_{1}\ra H_{0}$ is proper.
\end{dfn}

\Rem
Note that this implies that the other structure maps $\trg$, $\uni$ and
$\mlt$ of $H$ are also proper. Also note that if $H$ is proper then so is
the source-connected subgroupoid $H^{(0)}$ of $H$.

\begin{theo} \label{iascg1}
If $\fg$ and $\fh$ are integrable Lie algebroids and if
$\fh$ has a source-simply connected source-compact integral,
then any semi-direct product $\fg\ltimes\fh$ is integrable.
\end{theo}

\Proof
Let $G$ be a source-simply connected integral of $\fg$ and $H$
a source-simply connected source-compact integral of $\fh$.
Let $\nabla$ be an action of $\fg$ on $\fh$ along
$q:H_{0}\ra G_{0}$ over $R$.
For any $X\in\Gamma\fg$ let $w(X)$ be the multiplicative vector
field on $H$ with $\cL_{w(X)}=\nabla_{X}$ and $w(X)|_{H_{0}}=R(X)$
(Theorem \ref{dla5}). In particular, the map $w:\Gamma\fg\ra\fX^{\mlt}(H)$
is a $\eC(G_{0})$-linear homomorphism of Lie algebras, hence it
induces a bundle map $\bar{w}:(q\com\trg)^{\ast}\fg\ra T(H_{1})$
by $\bar{w}(h,u)=w(X)_{h}$, for any $X\in\Gamma\fg$ with
$X_{\pi(u)}=u$.

Now define a foliation $\cF$ of the manifold
$$ M=H_{1}\times_{G_{0}}G_{1}=\{(h,g)\in H_{1}\times G_{1}\,|\,
   q(\trg(h))=\trg(g)\} $$
by
$$ \cF_{(h,g)}=\{(Y_{\trg(h)}h + w(X)_{h},X_{\trg(g)}g)\,|\,
   X\in\Gamma\fg,\,Y\in\Gamma\fh\}\subset T_{(h,g)}M\;.$$
Note that $(Y_{\trg(h)}h + w(X)_{h},X_{\trg(g)}g)$
is indeed tangent to $M$ because
\begin{eqnarray*}
d(q\com\trg)(Y_{\trg(h)}h + w(X)_{h})
& = & dq(d\trg(w(X)_{h}))
  =   dq((R(X))_{\trg(h)})        \\
& = & \anchor(X)_{q(\trg(h))}
  =   \anchor(X)_{\trg(g)}
  =   d\trg(X_{\trg(g)}g)\;.
\end{eqnarray*}
The fact that $w$ is a $\eC(G_{0})$-linear homomorphism implies that
$\cF$ is a subbundle of $TM$, with $\dim(\cF)=\rank(\fg)+\rank(\fh)$.
We have to show that $\cF$ is involutive. For any
$X,X'\in\Gamma\fg$ and $Y,Y'\in\Gamma\fh$ we have
\begin{eqnarray}
&   & \!\!\!\!\!\!\!\!\!\!\!\!\!\!\!\!\!\!
      [(Y+w(X),X),(Y'+w(X'),X')] \nonumber \\
& = & ([Y+w(X),Y'+w(X')],[X,X']) \nonumber \\
& = & ([Y,Y']+[w(X),Y']+[Y,w(X')]+[w(X),w(X')],[X,X']) \nonumber \\
& = & ([Y,Y']+\nabla_{X}(Y')-\nabla_{X'}(Y)+w([X,X']),[X,X'])
      \label{eq2}\;,
\end{eqnarray}
and this is again a section of $\cF$ since
$[Y,Y']+\nabla_{X}(Y')-\nabla_{X'}(Y)\in\Gamma\fh$.
Here we denoted the invariant vector field on $G_{1}$ corresponding to
$X$ again by $X$, and the same for $X'$, $Y$ and $Y'$.

When we view the foliation $\cF$ of $M$ as a Lie algebroid over $M$,
there is a map of Lie algebroids
$$ \Phi:\cF\lra \fg\ltimes\fh $$
over the map $\trg\com\pr_{1}$, sending an $\cF$-tangent vector field
$(Y+ w(X),X)$ to $q^{\ast}X\oplus Y$. Indeed, the map $\Phi$ clearly
preserves the anchor, and it preserves the bracket
by Equation (\ref{eq2}). Note, in addition, that $\Phi$ restricts to an
isomorphism on each fiber. The algebroid $\cF$ is integrable by
$\Mon(M,\cF)$.

The manifold $M$ comes equipped with a natural right $H$-action along
the map $\src\com\pr_{1}$, which is principal with respect
to the projection $\trg\times\id:M\ra H_{0}\times_{G_{0}}G_{1}$.
Also, $M$ has a natural right
$G$-action along $\src\com\pr_{2}$, which is principal with
respect to the projection $\pr_{1}:M\ra H_{1}$.
These two actions commute with each other, and together they make
$M$ into a principal $(H\times G)$-bundle over $H_{0}$.
Note that the foliation $\cF$ is tangent to the fibers of
the map $\src\com\pr_{2}$, so the $G$-action lifts to
a $G$-action on the monodromy groupoid $\Mon(M,\cF)$.
On the other hand, the foliation $\cF$ is in general not tangent
to the fibers of $\src\com\pr_{1}$. Despite this we will
show that the $H$-action on $M$
may be lifted to an $H$-action on $\Mon(M,\cF)_{1}$.
However, this is not an action of $H$ on the groupoid in the sense
of Subsection \ref{sec:ag}.
We will show that $\Mon(M,\cF)$ can be factored by $G$ and $H$
with respect to these two actions to give a Lie groupoid $K$
over $H_{0}$ with the same $\src$-fibers as $\Mon(M,\cF)$.
It will then be clear from
the construction and the properties of the map $\Phi$
that $K$ integrates $\fg\ltimes\fh$.

To describe the $H$-action,
recall that source-connectedness of $H$ implies
that $q\com\src=q\com\trg$, so $H$ is a family of
Lie groupoids over $G_{0}$. We can take the pull-back of this family
along the target map $\trg:G_{1}\ra G_{0}$
to get a family of groupoids over $G_{1}$, as in the following diagram:
\begin{equation}
\xymatrix{
H_{1}\times_{H_{0}}H_{1}\times_{G_{0}}G_{1} \ar[r]^-{\mlt\times\id} &
H_{1}\times_{G_{0}}G_{1} \ar@<2pt>[r]^{\src\times\id} \ar@<-2pt>[r]_{\trg\times\id} &
H_{0}\times_{G_{0}}G_{1} \ar[r]^-{\pr_{2}} \ar@/_1pc/@<-4pt>[l]_{\uni\times\id} &
G_{1} \\
} \label{diag}
\end{equation}
The first pull-back here consists of $(h,h',g)\in H_{1}\times H_{1}\times G_{1}$
satisfying $\src(h)=\trg(h')$ and $\trg(h)=\trg(g)$.

Consider now the foliation $\cG$ on $G_{1}$ given by $\src$-fibers.
The action of $\fg$ on $H_{0}$ defines a foliation $\cG^{(0)}$ of
$H_{0}\times_{G_{0}}G_{1}$ given by
$$ \cG^{(0)}_{(y,g)}=\{(R(X)_{y},X_{\trg(g)}g)\,|\,X\in\Gamma\fg\}\;.$$
Similarly, define a foliation $\cG^{(1)}$ of
$H_{1}\times_{G_{0}}G_{1}$ by
$$ \cG^{(1)}_{(h,g)}=\{(w(X)_{h},X_{\trg(g)}g)\,|\,X\in\Gamma\fg\}\;,$$
and a foliation $\cG^{(2)}$ of
$H_{1}\times_{H_{0}}H_{1}\times_{G_{0}}G_{1}$ by
$$ \cG^{(2)}_{(h,h',g)}=\{(w(X)_{h},w(X)_{h'},X_{\trg(g)}g)\,|\,X\in\Gamma\fg\}\;.$$
Observe that, since each $w(X)$ is a multiplicative vector field on $H$, all the maps
in Diagram (\ref{diag}) map leaves to leaves.
Now note that the assumption that $H$ is source-compact implies that the maps
$\src\times\id$, $\trg\times\id$ and $\mlt\times\id$ of Diagram (\ref{diag})
are proper. By the argument given in the proof of Theorem \ref{iapm1}
it follows that each leaf of $\cG^{(1)}$ projects
along $\trg\times\id$ (and also along $\src\times\id$)
onto a leaf of $\cG^{(0)}$
as a covering projection. The same is true for the leaves
of $\cG^{(2)}$ with respect to the projection $\mlt\times\id$.
Note also that $\cG^{(1)}$ is a subfoliation of $\cF$.
There is another foliation $\cK$
of $H_{1}\times_{H_{0}}H_{1}\times_{G_{0}}G_{1}$ given by
$$ \cK_{(h,h',g)}=\{(w(X)_{h}+Y_{\trg(h)}h,w(X)_{h'},X_{\trg(g)}g)\,|\,
   X\in\Gamma\fg,\,Y\in\Gamma\fh\}\;.$$
The map $\mlt\times\id$ maps the leaves of $\cK$ onto the leaves of $\cF$
as a covering projection.

Now take any path $\gamma$ in a leaf of $\cF$ from $(h_{0},g_{0})$ to
$(h_{1},g_{1})$. Observe first that
$\bar{\gamma}=(\src\times\id)\com\gamma$ is a path tangent to
the foliation $\cG^{(0)}$. Hence for any $h\in H(\oo,\src(h_{0}))$ there
exists a unique lift $\tilde{\gamma}$ of $\bar{\gamma}$
along $\trg\times\id$
tangent to $\cG^{(1)}$ with $\tilde{\gamma}(0)=(h,g_{0})$.
Put $\delta=\pr_{1}\com\tilde{\gamma}$.
Now we may use the $H$-action on $M$ to define a new path
$\gamma h$ in $M$ by
$$ (\gamma h)(t)=\gamma(t) \delta(t)\;.$$
The fact that $\mlt\times\id$ maps $\cK$ to $\cF$
implies that $\gamma h$ is again a path in a leaf of $\cF$.
Note that $(\trg\times\id)\com (\gamma h)=(\trg\times\id)\com \gamma$.
It $h'$ is another arrow in
$H(\oo,\src(h))$, we lift $\bar{\gamma}$ to $\tilde{\gamma}$
along $\trg\times\id$ as before,
and we lift $\overline{\gamma h}$ along $\trg\times\id$
to $\tilde{\gamma}'$ tangent to $\cG^{(1)}$ with
$\tilde{\gamma}'(0)=(h',g_{0})$.
These give us
$\delta=\pr_{1}\com\tilde{\gamma}$ and $\delta'=\pr_{1}\com\tilde{\gamma}'$.
It follows that $(\delta,\delta',\pr_{2}\com\gamma)$ is
a path in a leaf of $\cG^{(2)}$, and its projection along
$\mlt\times\id$ is a path in a leaf of $\cG^{(1)}$ which lifts
$\bar{\gamma}$ with value $(hh',g_{0})$ at $t=0$.
Therefore
$$ \gamma (hh')=(\gamma h) h'\;.$$
A similar argument shows that
\begin{eqnarray}
(\gamma'\gamma)h=(\gamma' \tilde{h})(\gamma h)\;,
\label{eq4}
\end{eqnarray}
where $\tilde{h}$ is now the unique arrow in $H$ satisfying
$(\gamma h)(1)=\gamma(1)\tilde{h}$ (or $\tilde{h}=\delta(1)$ for
$\delta$ as above).

Finally, for any arrow $\sigma\in\Mon(M,\cF)_{1}$ we may define
$$ \sigma h \in \Mon(M,\cF)_{1} $$
by $\sigma h=[\gamma h]$, where $\gamma$ is any path representing
$\sigma$, i.e. $\sigma=[\gamma]$. This is well-defined since the
definition is given by path-lifting along a covering projection.
The properties mentioned above
imply that this defines an $H$-action on $\Mon(M,\cF)_{1}$ along the map
$\epsilon=\src\com\pr_{1}\com\src:\Mon(M,\cF)\ra H_{0}$.
This action commutes with the $G$-action, and we may take the
quotient $K=\Mon(M,\cF)/G/H$, which is a smooth manifold because
it can be identified with the pull-back of the source map of $\Mon(M,\cF)$
along $H_{0}\ra M$. Using Equation (\ref{eq4}) it is easy to check that
$K$ is a Lie groupoid over $H_{0}$.
Now $\Phi$ induces an isomorphism from the Lie algebroid of $K$ to
$\fg\ltimes\fh$ over $H_{0}$.
\eop
\vspace{3mm}

\Rem
Note that the assumption that $H$ is source-compact can be replaced by
the assumption that $H$ is Hausdorff and that all the vector fields
$w(X)$ are complete.

\section*{References}

\begin{list}{[\arabic{enumi}]}
{\usecounter{enumi}\settowidth\labelwidth{[99]}
\leftmargin\labelwidth\advance\leftmargin\labelsep}

\small

\bibitem{AM}
R. Almeida, P. Molino, Suites d'Atiyah et feuilletages transversalement complets.
{\em C. R. Acad. Sci. Paris 300,} 1985, pp. 13--15.

\bibitem{A}
M. F. Atiyah, Complex analytic connections in fibre bundles.
{\em Trans. Amer. Math. Soc. 85,} 1957, pp. 181--207.

\bibitem{BM}
R. Brown, O. Mucuk, The monodromy groupoid of a Lie groupoid.
{\em Cahiers Topologie Géom. Différentielle Catég. 36,} 1995, pp. 345--369.

\bibitem{CdSW}
A. Cannas da Silva, A. Weinstein, {\em Geometric Models for Noncommutative Geometry.}
Berkeley Mathematics Lecture Notes 10, American Mathematical Society, Providence, 1999.

\bibitem{CF}
A. S. Cattaneo, G. Felder,
Poisson sigma models and symplectic groupoids.
{\em arXiv:}\linebreak[0]~math.SG\linebreak[0]/0003023.

\bibitem{D}
P. Dazord, Groupo\"{\i}de d'holonomie et g\'{e}om\'{e}trie globale.
{\em C. R. Acad. Sci. Paris 324,} 1997, pp. 77--80.

\bibitem{Debord}
C. Debord, Groupo\"{\i}des d'holonomie de feuilletages singuliers.
{\em C. R. Acad. Sci. Paris 330,} 2000, pp. 361--364.

\bibitem{DL}
A. Douady, M. Lazard, Espaces fibr\'{e}s en alg\`{e}bres de Lie et en groupes.
{\em Invent. Math. 1,} 1966, pp. 133--151.

\bibitem{GHV}
W. Greub, S. Halperin, R. Vanstone, {\em Connections, Curvature and Cohomology.}
Pure and applied mathematics: a series of monographs and textbooks 47 Vol. I.,
Academic Press, New York, 1978.

\bibitem{HM}
P. J. Higgins, K. C. H. Mackenzie, Algebraic constructions in the category of Lie
algebroids. {\em J. Algebra 129,} 1990, pp. 194--230.

\bibitem{K}
I. M. Kontsevich, Deformation quantization of Poisson manifolds.
{\em arXiv:}\linebreak[0]~q-alg\linebreak[0]/9709040.

\bibitem{L}
N. P. Landsman, {\em Mathematical topics between classical and quantum mechanics.}
Springer Monographs in Mathematics, Springer-Verlag, New York, 1998

\bibitem{M}
K. C. H. Mackenzie, {\em Lie groupoids and Lie algebroids in Differential Geometry.}
London Mathematical Society Lecture Notes Series 124, Cambridge, 1987.

\bibitem{M1}
K. C. H. Mackenzie, Lie algebroids and Lie pseudoalgebras. {\em Bull. London Math.
Soc. 27,} 1995, pp. 97--147.

\bibitem{MX2}
K. C. H. Mackenzie, P. Xu, Lie bialgebroids and Poisson groupoids. {\em Duke
Math. J. 73,} 1994, pp. 415--452.

\bibitem{MX1}
K. C. H. Mackenzie, P. Xu, Classical lifting processes and multiplicative vector
fields. {\em Quart. J. Math. Oxford 49,} 1998, pp. 59--85.

\bibitem{MX}
K. C. H. Mackenzie, P. Xu, Integration of Lie bialgebroids. {\em Topology 39,}
2000, pp. 445--467.

\bibitem{Mo}
P. Molino, {\em Riemannian foliations.} Birkhauser, Boston, 1988.

\bibitem{N}
V. Nistor, Groupoids and the integration of Lie algebroids.
{\em arXiv:}\linebreak[0]~math.SG\linebreak[0]/0004084.

\bibitem{NWX}
V. Nistor, A. Weinstein, P. Xu,
Pseudodifferential operators on differential groupoids.
{\em Pacific J. Math. 189,} 1999, pp. 117--152.

\bibitem{NT}
R. Nest, B. Tsygan,
Deformations of symplectic Lie algebroids, deformations of holomorphic symplectic structures,
and index theorems. {\em arXiv:}\linebreak[0]~math.QA\linebreak[0]/9906020.

\bibitem{Pa}
R. Palais, {\em A Global Formulation of the Lie Theory of Transformation Groups.}
Memoirs of the American Mathematical Society 22, Providence, 1957.

\bibitem{Ph}
J. Phillips, The holonomic imperative and the homotopy groupoid of a foliated
manifold. {\em Rocky Mountain J. Math. 17,} 1987, pp. 151--165.

\bibitem{Pr}
J. Pradines, Th\'{e}orie de Lie pour les groupo\"{\i}des diff\'{e}rentiables. Calcul
diff\'{e}rentiel dans la cat\'{e}gorie des groupo\"{\i}des infinit\'{e}simaux.
{\em C. R. Acad. Sci. Paris 264,} 1967, pp. 245--248.

\bibitem{Pr1}
J. Pradines, Troisi\`{e}me th\'{e}or\`{e}me de Lie pour les groupo\"{\i}des
diff\'{e}rentiables. {\em C. R. Acad. Sci. Paris 267,} 1968, pp. 21--23.

\bibitem{S}
J.-P. Serre, {\em Lie Algebras and Lie Groups}. Springer-Verlag, Berlin Heidelberg 1992.

\bibitem{T}
F. Takens, Derivations of vector fields. {\em Compositio Math. 26,} 1973, pp. 151--158.

\bibitem{W}
H. Winkelnkemper, The graph of a foliation. {\em Ann. Global Anal. Geom. 1,} 1983,
pp. 51--75.

\normalsize

\end{list}

\end{document}